# FINITE-DIMENSIONAL APPROXIMATION FOR THE DIFFUSION COEFFICIENT IN THE SIMPLE EXCLUSION PROCESS

By Milton Jara

*IMPA*

We show that for the mean zero simple exclusion process in $\mathbb{Z}^d$ and for the asymmetric simple exclusion process in $\mathbb{Z}^d$ for $d \geq 3$, the self-diffusion coefficient of a tagged particle is stable when approximated by simple exclusion processes on large periodic lattices. The proof depends on a similar stability property of the Sobolev inner product associated with the operator.

**1. Introduction.** In [1], Kipnis and Varadhan proved an invariance principle for the position of a marked particle in a symmetric simple exclusion process in equilibrium. Their proof relies on a central limit theorem for additive functionals of a Markov process. Later, this result was generalized to the mean zero simple exclusion process (see [8]) and the asymmetric simple exclusion process in dimension $d \geq 3$ in [6].

The diffusion matrix of the limiting Brownian process is a function $D(\alpha)$ of the density of particles and is given by a variational formula.

The method of proof used by Kipnis and Varadhan works directly in infinite systems and naturally raises the issue of the stability of the diffusion coefficient under finite-dimensional approximations. More precisely, consider a finite-dimensional version of the simple exclusion process on the torus $\{-N, \ldots, 0, \ldots, N\}^d$. In order to obtain an ergodic process, fix the total number $K$ of particles. When $N$ is sufficiently large, the motion of a tagged particle on this finite system has a unique canonical lifting to $\mathbb{Z}^d$. In this manner, we obtain a process $X_N(t)$ with values in $\mathbb{Z}^d$. Let $D_{N,K}$ denote the variance of the limiting Brownian motion of the scaled process $\varepsilon X_N(t/\varepsilon^2)$ when $\varepsilon \to 0$. In this article, we prove that

$$\lim_{\substack{N \to \infty \\ K/(2N)^d \to \alpha}} D_{N,K} = D(\alpha)$$

for mean zero the simple exclusion process and for asymmetric the simple exclusion process in dimension $d \geq 3$.

This limit was first considered in [3] for symmetric the simple exclusion process. The proof presented there follows from a variational formula for the diffusion coefficient that depends on the Sobolev dual norm associated with the generator of the process and from a convergence result for the Sobolev dual norms of the finite-dimensional approximations. Let $h, g$ be local functions with mean zero with respect to all the Bernoulli product measures $\mu_\alpha$ that assign density $\alpha$ to each coordinate. Denote by $\langle \cdot, \cdot \rangle_\alpha$ the inner product on $\mathcal{L}^2(\mu_\alpha)$. Let $\mu_{N,K}$ be the uniform measure over the configurations with $K$ particles on the torus $\{-N, \ldots, 0, \ldots, N\}^d$ and $\langle \cdot, \cdot \rangle_{\mu_{N,K}}$ the inner product on $\mathcal{L}^2(\mu_{N,K})$. Let $L$ (resp. $L_N$) be the generator of the process in $\mathbb{Z}^d$ (resp. the torus). Suppose for a moment that $(-L)^{-1}g$ exists and is local. Then

$$\lim_{\substack{N \to \infty \\ K/(2N)^d \to \alpha}} \langle h, (-L_N)^{-1} g \rangle_{N,K} = \langle h, (-L)^{-1} g \rangle_\alpha$$

because of the equivalence of ensembles and the fact that $(-L)^{-1}g$ is local. The desired result will be a consequence of a generalization of this result for a larger class of functions $h$, $g$.

The rest of this article is organized as follows. In Section 2, we present the definition and basic properties of the simple exclusion process. In Section 3, we introduce the Sobolev spaces associated with the process and prove Theorem 2, a general perturbative result about the convergence of finite approximations of a positive operator in Hilbert spaces. In Section 4, we prove the stability of the diffusion coefficient for the tagged particle using the Sobolev space techniques developed in Section 3. Finally, in Section 5, we check the hypothesis of Theorem 2 for the simple exclusion process.

**2. Notation and results.** Consider a probability measure $p(\cdot)$ of finite range on $\mathbb{Z}^d$, that is, $p(z) = 0$ if $|z|$ is sufficiently large. Suppose that $p(0) = 0$ and that the random walk with transition rate $p(\cdot)$ is irreducible, that is, the (finite) set $\{z; p(z) > 0\}$ generates the group $\mathbb{Z}^d$. The simple exclusion process associated with $p(\cdot)$ corresponds to the Markov process defined on $\mathcal{X} = \{0,1\}^{\mathbb{Z}^d}$, whose generator $\mathcal{L}_0$ acting on local functions $f$ is given by

$$\mathcal{L}_0 f(\eta) = \sum_{x,y \in \mathbb{Z}^d} p(y-x)\eta(x)(1-\eta(y))[f(\sigma^{xy}\eta) - f(\eta)].$$

Here, $\eta \in \mathcal{X}$ denotes a configuration of particles in $\mathbb{Z}^d$. In particular, $\eta(x) = 1$ if there is a particle at the site $x$ and $\eta(x) = 0$ otherwise, and $\sigma^{xy}\eta$ is the configuration obtained from $\eta$ by exchanging the occupation

numbers at $x$ and $y$:

$$\sigma^{xy}\eta(z) = \begin{cases} \eta(y), & \text{if } z = x, \\ \eta(x), & \text{if } z = y, \\ \eta(z), & \text{otherwise.} \end{cases}$$

If $p(z) = p(-z)$ for all $z$, then the process will be called *symmetric*; if $\sum zp(z) = 0$, it will be called *of mean zero* and if $\sum_{z \in \mathbb{Z}^d} zp(z) = m \neq 0$, the process will be called *asymmetric*.

For each $\alpha \in [0,1]$, let $\nu_\alpha$ denote the Bernoulli product measure in $\mathcal{X}$, that is, the product measure such that $\nu_\alpha[\eta(x) = 1] = \alpha$ for each $x \in \mathbb{Z}^d$. It is not hard to prove that $\nu_\alpha$ is an invariant measure for the process generated by $\mathcal{L}_0$.

In this model, particles are indistinguishable. In order to study the time evolution of a single particle, we proceed in the following way. Let $\eta \in \mathcal{X}$ be an initial state with a particle at the origin [i.e., $\eta(0) = 1$]. Tag this particle and let $\eta_t$ (resp. $X_t$) be the time evolution of the exclusion process starting from $\eta$ (resp. the tagged particle starting from $x = 0$). Let $\xi_t(x) = \eta_t(x + X_t)$ be the process as seen by the tagged particle. We call $\xi_t$ the *environment process*.

It is clear that $X_t$ is not a Markov process due to the interaction between the tagged particle and the environment, but $(\eta_t, X_t)$ and $\xi_t$ are Markov processes, the latter being defined in the state space $\mathcal{X}_* = \{0,1\}^{\mathbb{Z}^d_*}$, where $\mathbb{Z}^d_* = \mathbb{Z}^d \setminus \{0\}$. The generator of the process $\xi_t$, acting on local functions $f$, is given by $L = L_0 + L_\tau$, where

$$L_0 f(\xi) = \sum_{x,y \in \mathbb{Z}^d_*} p(y-x)\xi(x)(1-\xi(y))[f(\sigma^{xy}\xi) - f(\xi)],$$

$$L_\tau f(\xi) = \sum_{z \in \mathbb{Z}^d_*} p(z)(1-\xi(z))[f(\tau_z\xi) - f(\xi)].$$

The first part of the generator, $L_0$, takes into account the jumps of the environment (i.e., all particles but the tagged one), while the second part takes into account the jumps of the tagged particle.

In this formula, $\tau_z\xi$ is the configuration obtained by making the tagged particle (at the origin) jump to site $z$ and then bringing it back to the origin with a translation:

$$\tau_z\xi(x) = \begin{cases} 0, & \text{if } x = -z, \\ \xi(x+z), & \text{if } x \neq -z. \end{cases}$$

For the process $\xi_t$, we have a one-parameter family of invariant ergodic measures $\{\mu_\alpha\}_{\alpha \in [0,1]}$, where $\mu_\alpha$ is the Bernoulli product measure of density $\alpha$ defined on $\mathcal{X}_*$: $\mu_\alpha[\xi(x) = 1] = \alpha$ for all $x \in \mathbb{Z}^d_*$, independently for each site (see [5]).

Note that the position of the tagged particle can be calculated in terms of jump processes associated with $\xi_t$. Define $N_t^z$ as the number of translations by $z$ of $\xi_t$, that is, $N_t^z = N_{t-}^z + 1 \iff \xi_t = \tau_z \xi_{t-}$. Then $X_t = \sum_z z N_t^z$.

In this context, Kipnis and Varadhan proved a central limit theorem for the position of the tagged particle when the environment process is in equilibrium with distribution $\mu_\alpha$. They proved that $\varepsilon X_{t/\varepsilon^2}$ converges, when $\varepsilon$ goes to zero, to a Brownian motion with diffusion coefficient $D(\alpha)$, which can be described in terms of the Sobolev norms associated with the operator $L$ in $\mathcal{L}(\mu_\alpha)$.

This result has been generalized by Varadhan to the mean zero case (in any dimension) and by Sethuraman, Varadhan and Yau for the asymmetric case in dimension $d \geq 3$, in which case it is proved that $\varepsilon[X_{t/\varepsilon^2} - mt(1-\alpha)/\varepsilon^2]$ converges to a Brownian motion with diffusion coefficient $D(\alpha)$, given by

$$a^t D(\alpha) a = (1-\alpha) \sum_{z \in \mathbb{Z}_*^d} (z \cdot a)^2 p(z) - 2 \langle w_a, (-L)^{-1} v_a \rangle_\alpha, \tag{1}$$

where $a \in \mathbb{R}^d$, $\langle \cdot, \cdot \rangle_\alpha$ is the inner product on $\mathcal{L}^2(\mu_\alpha)$ and $v_a$, $w_a$ are local functions defined by

$$v_a = \sum_{z \in \mathbb{Z}^d} (z \cdot a) p(z) [\alpha - \eta(z)],$$

$$w_a = \sum_{z \in \mathbb{Z}^d} (z \cdot a) p(z) [\alpha - \eta(-z)].$$

In general, $L$ is not an invertible operator and so the meaning of (1) must be clarified. This will be done in Sections 3 and 4.

Let $N$ be a positive integer and define $T_N^d = \{-N, \ldots, 0, \ldots, N\}^d$, the $d$-dimensional discrete torus of $(2N)^d$ points, with $-N$ and $N$ identified. Using the same probability measure $p(\cdot)$, we can define a simple exclusion process evolving in $T_N^d$. The space state will now be $\mathcal{X}_N = \{0,1\}^{T_N^d}$ and the generator $\mathcal{L}_N$ acting on any function $f$ will be given by

$$\mathcal{L}_N f(\xi) = \sum_{x,y \in T_N^d} p(y-x) \eta(x)(1-\eta(y))[f(\sigma^{xy}\eta) - f(\eta)].$$

In the same way, it is possible to define the environment process in the torus $T_{N,*}^d = T_N^d \setminus \{0\}$. In this case, the environment, as seen by the tagged particle, is a Markov process evolving in the space $\mathcal{X}_{N,*} = \{0,1\}^{T_{N,*}^d}$ and generated by the operator $L_N = L_{0,N} + L_{\tau,N}$, where

$$L_{0,N} f(\xi) = \sum_{x,y \in T_{N,*}^d} p(y-x) \xi(x)(1-\xi(y))[f(\sigma^{xy}\xi) - f(\xi)],$$

$$L_{\tau,N} f(\xi) = \sum_{z \in T_{N,*}^d} p(z)(1 - \xi(z))[f(\tau_z \xi) - f(\xi)].$$

It is clear, by the conservation of the number of particles, that for $0 < K \leq (2N)^d$, the probability measure $\mu_{N,K}$, uniform over the set $\mathcal{X}_{N,K} = \{\xi \in \mathcal{X}_{N,*}; \sum_{x \in T_{N,*}^d} \xi(x) = K - 1\}$ of configurations with $K$ particles, is an invariant ergodic measure for the process generated by $L_N$.

For $N$ sufficiently large, it is possible to lift the motion of the tagged particle to $\mathbb{Z}^d$. Let $X_t^N$ denote the position of the tagged particle in $\mathbb{Z}^d$. It is not hard to prove an invariance principle for $X_t^N : \varepsilon[x_{t/\varepsilon^2}^N - mt(1 - \alpha_{N,K})/\varepsilon^2]$ converges to a Brownian motion of variance $D_{N,K}$, given by

$$\begin{aligned} a^t D_{N,K} a = (1 - \alpha_{N,K}) \sum_{z \in \mathbb{Z}^d} (a \cdot z)^2 p(z) \\ - 2 \langle w_a - \langle w_a \rangle_{N,K}, L_N^{-1}(v_a - \langle v_a \rangle_{N,K}) \rangle_{N,K}. \end{aligned}$$
(2)

In this formula, $\langle \cdot, \cdot \rangle_{N,K}$ (resp. $\langle \cdot \rangle_{N,K}$) stands for the inner product on $\mathcal{L}(\mu_{N,K})$ (resp. the mean with respect to $\mu_{N,K}$) and

$$\alpha_{N,K} = \frac{K-1}{(2N)^d - 1}.$$

Note that for $f : \mathcal{X}_{N,K} \to \mathbb{R}$ with $\langle f \rangle_{N,K} = 0$, $L_N^{-1} f$ is well defined. In fact, for $f$, we have

$$\begin{aligned} \mathcal{D}_{N,K}(f) &= \langle f, -L_N f \rangle_{N,K} \\ &= \tfrac{1}{4} \sum_{x,y \in T_{N,*}^d} (p(y-x) + p(x-y)) \int [f(\sigma_{x,y}\eta) - f(\eta)]^2 \, d\mu_{N,K}. \end{aligned}$$

In particular, $L_N f = 0$ if and only if $f$ is constant. Also, $L_N$ is an invertible operator in $\mathcal{C}_{0,N,K} = \{f; \langle f \rangle_{N,K} = 0\}$.

For the symmetric simple exclusion process, Landim, Olla and Varadhan [3] proved that $D_{N,K} \to D(\alpha)$ if $\alpha_{N,K} \to \alpha$. In this article, we extend this result to the asymmetric case, as given by the following theorem:

THEOREM 1. *For the mean zero simple exclusion process (in any dimension) and for the asymmetric simple exclusion process in dimension $d \geq 3$, $D_{N,K} \to D(\alpha)$ if $\alpha_{N,K} \to \alpha$.*

**3. The Sobolev spaces $H_1$, $H_{-1}$.** In this section, we prove the stability of the $H_{-1}$ norm under finite approximations. We discuss it in the more general context of functional analysis because it is a general result that can be applied to many models of interacting particle systems and we will be used repeatedly in the sequel.

Let $H$ be a real Hilbert space with inner product $\langle \cdot, \cdot \rangle$. An operator (not necessarily bounded) $L: D(L) \subseteq H \to H$ is called *positive* if $\langle g, Lg \rangle > 0$ for all $g \in D(L) \setminus \{0\}$.

Given a positive closed operator $L$, we define, for $f \in D(L)$,

$$\|f\|_1^2 =: \langle f, Lf \rangle.$$

It is easy to see that $\|\cdot\|_1$ defines a norm on $D(L)$ that satisfies the parallelogram rule. Therefore, $\|\cdot\|_1$ can be extended to an inner product on $D(L)$. Define $\mathcal{H}_1 = \mathcal{H}_1(L)$, the Sobolev space associated with the operator $L$, as the completion of $D(L)$ under $\|\cdot\|_1$.

In the same way, we see that

$$\|g\|_{-1}^2 =: \sup_{f \in D(L)} \{2\langle g, f \rangle - \langle f, Lf \rangle\}$$

defines a norm on the set $\{g \in H; \|g\|_{-1} < \infty\}$ that can be extended to an inner product. Define $\mathcal{H}_{-1}$ as the completion of this set under $\|\cdot\|_{-1}$.

In the next proposition, some well-known properties of the spaces $\mathcal{H}_1, \mathcal{H}_{-1}$ are listed.

PROPOSITION 1. *For $f \in H \cap \mathcal{H}_1$, $g \in H \cap \mathcal{H}_{-1}$, we have*

(i) $\|g\|_{-1} = \sup_{h \in D(L) \setminus \{0\}} \frac{\langle h, g \rangle}{\|h\|_1}$,
(ii) $|\langle f, g \rangle| \leq \|f\|_1 \|g\|_{-1}$,
(iii) $\|f\|_1 \leq \|Lf\|_{-1}$.

PROOF. For (i),

$$\sup_{h \in D(L)} \{2\langle g, h \rangle - \langle h, Lh \rangle\} = \sup_{\|h\|_1 = 1} \sup_{\alpha \in \mathbb{R}} \{2\alpha \langle g, h \rangle - \alpha^2\}$$

$$= \sup_{\|h\|_1 = 1} \langle g, h \rangle^2.$$

For (ii),

$$\|g\|_{-1} = \sup_{h \in D(L) \setminus \{0\}} \frac{|\langle g, h \rangle|}{\|h\|_1} \geq \frac{|\langle g, f \rangle|}{\|f\|_1}.$$

For (iii),

$$\|Lf\|_{-1}^2 = \sup_{h \in D(L)} \{2\langle Lf, h \rangle - \langle h, Lh \rangle\} \geq \langle f, Lf \rangle. \qquad \square$$

From property (i), it can be concluded that $\mathcal{H}_{-1}$ is the dual of $\mathcal{H}_1$ with respect to $H$. Thanks to property (ii), the inner product $\langle \cdot, \cdot \rangle$ can be extended to a continuous bilinear form $\langle \cdot, \cdot \rangle : \mathcal{H}_{-1} \times \mathcal{H}_1 \to \mathbb{R}$. Property (iii) ensures that the operator $L^{-1}: \text{Im}(L) \cap \mathcal{H}_{-1} \to \mathcal{H}_1$ is bounded, from which it can be continuously extended to an operator defined on the closure of $\operatorname{Im}(L) \cap \mathcal{H}_{-1}$ under $\|\cdot\|_{-1}$.

If the operator $L$ is symmetric, that is, if $\langle f, Lg \rangle = \langle Lf, g \rangle$ for $f, g \in D(L)$, then the inequality in (iii) becomes an equality and $L$ can be extended to an isometry from $\mathcal{H}_1$ to $\mathcal{H}_{-1}$ (not necessarily surjective).

Let $\{H_n\}_n$ be an increasing sequence of finite-dimensional subspaces of $H$ and define $\operatorname{Loc} = \operatorname{Loc}(H) =: \bigcup_n H_n$. Suppose that Loc is a kernel for $L$, that is, the closure of the operator $L$ restricted to Loc is the operator $L$ itself. Suppose also that Loc is a kernel for the adjoint $L^*$ of $L$. Consider on each subspace $H_n$, an inner product $\langle \cdot, \cdot \rangle_n$ such that for all $f, g \in \operatorname{Loc}$,

$$\lim_{n \to \infty} \langle f, g \rangle_n = \langle f, g \rangle,$$

where $\langle f, g \rangle_n$ is well defined for $n$ sufficiently large.

A sequence $\{L_n\}_n$ of operators is called a *finite approximation* of $L$ if:

(i) $L_n : H_n \to H_n$;
(ii) $\langle f, L_n f \rangle_n > 0$ for $f \in H_n \setminus \{0\}$;
(iii) for all $f \in \operatorname{Loc}$, there exist $n_0 \in \mathbb{N}$ such that $L_n f = Lf$ for $n \geq n_0$;
(iv) if $L$ is a symmetric operator, then $L_n$ is also a symmetric operator.

On $H_n$, define the $\|\cdot\|_{1,n}$, $\|\cdot\|_{-1,n}$ norms associated with $L_n$ as before:

$$\|f\|_{1,n}^2 = \langle f, L_n f \rangle_n,$$
$$\|f\|_{-1,n}^2 = \sup_{g \in H_n} \{2 \langle f, g \rangle_n - \langle g, L_n g \rangle_n\}.$$

Observe that $\operatorname{Ker}(L_n) = \{0\}$, from which it follows that $L_n$ is invertible. The purpose of this section is to establish sufficient conditions to ensure that

$$(\star) \qquad \lim_{n \to \infty} \langle h', L_n^{-1} h \rangle_n = \langle h', L^{-1} h \rangle$$

for $h, h' \in \operatorname{Loc} \cap \mathcal{H}_{-1}$ with $h$ in the closure of $\operatorname{Im}(L) \cap \mathcal{H}_{-1}$.

While $L_n^{-1} h$ is always well defined, $h$ might not be in the image of $L$, in which case the left-hand side of this equality would not be well defined. However, when $h$ is in the closure of $\operatorname{Im}(L) \cap \mathcal{H}_{-1}$ under $\|\cdot\|_{-1}$, the product $\langle h', L^{-1} h \rangle$ can be defined by continuity. Recall that the product $\langle h', L_n^{-1} h \rangle_n$ is well defined for $n$ sufficiently large, because $h, h' \in \operatorname{Loc}$; each time a limit like the one appearing in $(\star)$ is considered, this fact must be taken into account.

The next theorem is a perturbative result asserting that if $(\star)$ is satisfied for an operator $S_0$ (and a suitable finite approximation $\{S_{0,n}\}_n$ of $S_0$), then it is also satisfied for a class of perturbations of $S_0$.

THEOREM 2. *Let $L$ be a positive closed operator. Let $S_0 : D(S_0) \subseteq H \to H$ be a symmetric positive operator such that* Loc *is a kernel for $S_0$ and $\langle g, S_0 g \rangle \leq \langle g, Lg \rangle$. Let $\{S_{0,n}\}_n$ be a finite approximation of $S_0$ such that $\langle f, S_{0,n} f \rangle_n \leq \langle f, L_n f \rangle_n$ for all $f \in H_n$. Define the norms $\|\cdot\|_{0,1}$, $\|\cdot\|_{0,-1}$ (resp. $\|\cdot\|_{0,1,n}, \|\cdot\|_{0,-1,n}$) associated to $S_0$ (resp. $S_{0,n}$) as before. Consider $h, h' \in \mathrm{Loc} \cap \mathcal{H}_{-1}$, with $h$ in the closure of $\mathrm{Im}(L) \cap \mathcal{H}_{-1}$.*

*Assume that*

(A) *For each $\varepsilon > 0$, there exists $g_\varepsilon \in \mathrm{Loc}$ such that*
$$\|h - Lg_\varepsilon\|_{0,-1} < \varepsilon.$$

(B)
$$\lim_{n \to \infty} \|h'\|_{0,-1,n} = \|h'\|_{0,-1}$$

*and for $u_\varepsilon = h - Lg_\varepsilon$,*
$$\lim_{n \to \infty} \|u_\varepsilon\|_{0,-1,n} = \|u_\varepsilon\|_{0,-1}.$$

*Then*
$$\lim_{n \to \infty} \langle h', L_n^{-1} h \rangle_n = \langle h', L^{-1} h \rangle.$$

PROOF. First, we observe that the operator $S_0$ (resp. $S_{0,n}$) is dominated by $L_0$ (resp. $L_{0,n}$), from which we have, for all $f$, the inequalities
$$\|f\|_{-1} \leq \|f\|_{0,-1},$$
$$\|f\|_{0,1} \leq \|f\|_1,$$
$$\|f\|_{-1,n} \leq \|f\|_{0,-1,n},$$
$$\|f\|_{0,1,n} \leq \|f\|_{1,n}.$$

Fix $\varepsilon > 0$ and let $u_\varepsilon = h - Lg_\varepsilon$ be chosen according to assumption (A). Then $Lg_\varepsilon = L_n g_\varepsilon$ for $n$ sufficiently large, from which it follows that $Lg_\varepsilon$ belongs to $H_n$ and
$$\langle h', L_n^{-1} h \rangle_n = \langle h', L_n^{-1}(u_\varepsilon + Lg_\varepsilon) \rangle_n$$
$$= \langle h', g_\varepsilon \rangle_n + \langle h', L_n^{-1} u_\varepsilon \rangle_n.$$

Since $h'$ and $g_\varepsilon$ are in Loc, we have
$$\langle h', g_\varepsilon \rangle_n \xrightarrow[n \to \infty]{} \langle h', g_\varepsilon \rangle.$$

We also have that
$$|\langle h', L_n^{-1} u_\varepsilon \rangle_n| \leq \|h'\|_{0,-1,n} \cdot \|L_n^{-1} u_\varepsilon\|_{0,1,n}$$
$$\leq \|h'\|_{0,-1,n} \cdot \|L_n^{-1} u_\varepsilon\|_{1,n}$$
$$\leq \|h'\|_{0,-1,n} \cdot \|u_\varepsilon\|_{-1,n}$$
$$\leq \|h'\|_{0,-1,n} \cdot \|u_\varepsilon\|_{0,-1,n}.$$

Therefore,
$$\limsup_{n\to\infty} |\langle h', L_n^{-1} u_\varepsilon \rangle_n| \leq \|h'\|_{0,-1} \cdot \|u_\varepsilon\|_{0,-1} \leq \varepsilon \cdot \|h'\|_{0,-1}.$$

On the other hand, $\langle h', L^{-1}h \rangle = \langle h', L^{-1}u_\varepsilon \rangle + \langle h', g_\varepsilon \rangle$ and
$$|\langle h', L^{-1}u_\varepsilon \rangle| \leq \|h'\|_{0,-1} \cdot \|L^{-1}u_\varepsilon\|_{0,1} \leq \varepsilon \cdot \|h'\|_{0,-1}.$$

Consequently,
$$\limsup_{n\to\infty} |\langle h', L^{-1}h \rangle - \langle h', L_n^{-1}h \rangle_n| \leq 2\varepsilon \|h'\|_{0,-1}. \qquad \square$$

**4. Proof of Theorem 1.** This section is organized as follows. First, we show in which sense the sequence $\{L_N\}_N$ is a finite approximation of the operator $L$. Once this has been done, the proof of Theorem 1 is reduced to the verification of the hypothesis of Theorem 2, as we will see. We then verify these hypotheses separately for the symmetric, mean zero and asymmetric simple exclusion process.

4.1. *Finite approximations for $L$.* Let $\alpha \in [0,1]$ be fixed. Let $\{K_N\}_N$ be a sequence such that as $N$ goes to infinity, $\alpha_{N,K_N} \to \alpha$, $K_N \to \infty$ and $(2N)^d - K_N \to \infty$. Hereafter, we omit the index $K_N$ if there is no risk of confusion. Let $f, g$ be in $\mathcal{L}^2(\mu_\alpha)$. First, we deal with irrelevant constants. We say that $f \sim g$ if $\int (f-g) \, d\mu_\alpha = 0$. Define $H = \mathcal{L}^2(\mu_\alpha)/\sim$. It is easy to see that $H$ is isomorphic to the set of functions with mean zero in $\mathcal{L}^2(\mu_\alpha)$. Let $\mathrm{Loc} = \mathrm{Loc}(H)$ be the set of local functions in $H$. We define $H_N \cong \mathcal{C}_{0,N,K_N}$ as follows. Consider the canonical projection $\pi_N : \mathcal{X}_* \to \mathcal{X}_{N,*}$. For $f \in \mathcal{C}_{0,N,K_N}$, define $\pi_N^{-1} f \in \mathrm{Loc}$ by
$$\pi_N^{-1} f(\eta) = \begin{cases} f(\pi_N \eta), & \text{if } \pi_N \eta \in \mathcal{X}_{N,K_N}, \\ 0, & \text{if } \pi_N \eta \notin \mathcal{X}_{N,K_N}. \end{cases}$$

Then $H_N = \pi_N^{-1}(\mathcal{C}_{0,N,K_N})/\sim$. It is not hard to see that $\mathrm{Loc} = \bigcup_N H_N$. In fact, for a local function $f$, denote by $\mathrm{supp}(f)$ the support of $f$. Then if $\mathrm{supp}(f) \subseteq T_{N,*}^d$, $\#\mathrm{supp}(f) < \min\{K_N, (2N)^d - K_N\}$, it follows that $f \in H_N$ and so, clearly, $H_N \subseteq \mathrm{Loc}$. On $H_N$, we define the inner product $\langle \cdot, \cdot \rangle_N$ induced by the measure $\mu_{N,K_N}$.

It is clear that for $f, g \in \mathrm{Loc}$ and $N$ sufficiently large (note that $f, g$ are not necessarily in $\mathcal{C}_{0,N,K_N}$),
$$\langle f, g \rangle_N = \int \left( f - \int f \, d\mu_{N,K_N} \right) \left( g - \int g \, d\mu_{N,K_N} \right) d\mu_{N,K_N}$$
$$= \int fg \, d\mu_{N,K_N} - \int f \, d\mu_{N,K_N} \int g \, d\mu_{N,K_N}.$$

We have already seen that the operator $-L_N$ is positive and it is clear that $-L_N f = Lf$ for $f \in \mathrm{Loc}$ and $N$ sufficiently large. From the ergodicity of $\mu_\alpha$ with respect to the process generated by $L$ and the fact that $L$ is a generator of a Markov process, we deduce that $\mathcal{D}_\alpha(f) = \langle f, -Lf \rangle_\alpha > 0$ if $f \neq 0$, from which we see that $-L$ is a positive operator. Consequently, $\{L_N\}_N$ would be a finite approximation of $-L$ if it were not for the fact that $H_N \subsetneq H_{N+1}$ [because $K_N$, $(2N)^d - K_N$ are not necessarily increasing sequences]. However, it is true that $H_N \subseteq H_M$ for $M$ sufficiently large, where $M$ depends both on the range of the transition probability $p(\cdot)$ and on the sequence $K_N$ [here, we use the fact that $K_N \to \infty$ and $(2N)^d - K_N \to \infty$]. Of course, Theorem 1 applies in this situation by taking subsequences or slightly modifying it to fit this case. Anyway, we will say that $\{L_N\}_N$ is a finite approximation of $L$.

Note that the inner product $\langle \cdot, \cdot \rangle_N$ is exactly the product appearing in equation (2). Comparing equations (1) and (2), it is clear that Theorem 1 follows from Theorem 2 applied to the operators $-L$ and $-L_N$. So, it only remains to find suitable operators $S_0$ and $\{S_{0,N}\}$, to compare with $L$ and $\{L_N\}$ and to check the hypothesis of Theorem 2 for them.

4.2. *Symmetric case.* Suppose that the transition probability $p(\cdot)$ is symmetric, that is, $p(x) = p(-x)$ for all $x \in \mathbb{Z}^d$. This case has been considered in [3], but in order to make the exposition clear, we here outline the proof in our setting.

Choose $S_0 = -L_0$ and $S_{0,N} = -L_{0,N}$, the part of the generator corresponding to jumps in the environment. It is clear that $\{S_{0,N}\}$ is a finite approximation of $S_0$ and that $\langle f, S_0 f \rangle_\alpha \leq \langle f, -Lf \rangle_\alpha$, $\langle g, S_{0,N} g \rangle_N \leq \langle g, -L_N g \rangle_N$. Conditions (A) and (B) of Theorem 2 are consequences in this case, of the next results, with we state as lemmas.

LEMMA 1. $w_a, v_a \in H_{0,-1}$ *and for all* $g \in \text{Loc}$, $Lg \in H_{0,-1}$.

PROOF. Following a criterion of Sethuraman and Xu [7], a sufficient condition for a local function $v$ to be in $H_{0,-1}$ is that $\langle v \rangle_\alpha = 0$ for all $\alpha \in [0,1]$. Therefore, it is enough to observe that for all $\alpha \in [0,1]$, $\langle w_a \rangle_\alpha = \langle v_a \rangle_\alpha = \langle Lg \rangle_\alpha = 0$. □

LEMMA 2. *If* $g \in \text{Loc}$ *and* $\langle g \rangle_\alpha = 0$ *for all* $\alpha \in [0,1]$, *then*
$$\lim_{N \to \infty} \|g\|_{0,-1,N} = \|g\|_{0,-1}.$$

PROOF. This is just a consequence of Corollaries 2.2 and 2.4 of [3] that are based on the so-called Liouville-D property of the lattice $\mathbb{Z}_*^d$. □

The following lemma is just Theorem 4.2 of [3]:

LEMMA 3. *If $v \in \text{Loc}$ and $\langle v \rangle_\alpha = 0$ for all $\alpha \in [0,1]$, then for all $\varepsilon > 0$, there exists $g_\varepsilon \in \text{Loc}$ such that*

$$\|v - Lg_\varepsilon\|_{0,-1} < \varepsilon.$$

Once these three lemmas are stated, by Theorem 2 we have the following result:

THEOREM 3. *For all $v \in \text{Loc}$ such that $\langle v \rangle_\alpha = 0$ for all $\alpha \in [0,1]$, we have*

$$\lim_{N \to \infty} \|v\|_{-1,N} = \|v\|_{-1}.$$

4.3. *Mean zero case.* Now, suppose that the transition probability has mean zero, that is, $\sum_z z p(z) = 0$. Define $S = -(L + L^*)/2$, $S_N = -(L_N + L_N^*)/2$, the symmetric part of the generator. A simple computation shows that

$$Sf(\xi) = \sum_{x,y \in \mathbb{Z}^d_*} s(y-x)\xi(x)(1-\xi(y))[f(\sigma^{xy}\xi) - f(\xi)]$$

$$+ \sum_{z \in \mathbb{Z}^d_*} s(z)(1-\xi(z))[f(\tau_z\xi) - f(\xi)]$$

and

$$S_N f(\xi) = \sum_{x,y \in T^d_{N,*}} s(y-x)\xi(x)(1-\xi(y))[f(\sigma^{xy}\xi) - f(\xi)]$$

$$+ \sum_{z \in T^d_{N,*}} s(z)(1-\xi(z))[f(\tau_z\xi) - f(\xi)],$$

where $s(x) = (p(x) + p(-x))/2$, the symmetrization of $p(\cdot)$. It is clear that $s(\cdot)$ is a symmetric, finite-range, irreducible transition probability, from which $S$ (resp. $S_N$) is the generator of a symmetric exclusion process in $\mathbb{Z}^d_*$ (resp. $T^d_{N,*}$). We choose $S_0 = S$ and $S_{0,N} = S_N$. As in the symmetric case, $S_{0,N_N}$ is a finite approximation of $S_0$ and, by definition, $\langle f, S_0 f \rangle = \langle f, -Lf \rangle$ and $\langle f, S_{0,N} f \rangle_N = \langle f, -L_N f \rangle_N$. Observe that in this case, $S_0$ and $-L$ generate the same Sobolev norms.

As in the symmetric case, we need to verify Assumptions (A) and (B) of Theorem 2. First, we need to prove that $w_a, v_a \in H_{-1}$ and for $g \in \text{Loc}$, that $Lg \in H_{-1}$. But this is true because $\langle v_a \rangle_\alpha = \langle w_a \rangle_\alpha = \langle Lg \rangle_\alpha = 0$ for all $\alpha \in [0,1]$, $H_{-1} \subseteq H_{0,-1}$ (in the notation of the previous subsection) and by the criterion of [7], $v_a, w_a, Lg \in H_{0,-1}$.

Assumption (B) of Theorem 2 then follows from Theorem 3. Therefore, in order to apply Theorem 2 to prove Theorem 1, it only remains to prove assumption (A). We state it as the following lemma:

LEMMA 4. *For all $v \in \mathrm{Loc}$ such that $\langle v \rangle_\alpha = 0$ for all $\alpha \in [0,1]$, and for all $\varepsilon > 0$, there exists $g_\varepsilon \in \mathrm{Loc}$ such that*

$$\|v - Lg_\varepsilon\|_{-1} < \varepsilon.$$

PROOF. In [8], Varadhan proved a sector condition for the mean zero exclusion process, which roughly states that the asymmetric part of the operator can be bounded by the symmetric part. More precisely, there exists a constant $C = C(p(\cdot))$ such that for all $f, g \in \mathrm{Loc}$,

$$\langle f, Lg \rangle_\alpha^2 \leq C \langle f, -Lf \rangle_\alpha \langle g, -Lg \rangle_\alpha.$$

In particular, $\|Lg\|_{-1}^2 \leq C\|g\|_1^2$, from which it follows that $L$ is a bounded and densely defined operator from $H_1$ to $H_{-1}$. So, it is enough to prove that $v \in L(H_1)$. To this end, we use the resolvent method. Let $h$ be in $H_{-1} \cap \mathrm{Loc}$. For each $\lambda > 0$, let $u_\lambda$ be the solution of the resolvent equation

$$\lambda u_\lambda - L u_\lambda = h.$$

This is always possible because $L$ is a negative operator in $\mathcal{L}^2(\mu_\alpha)$ and $u_\lambda \in D(L)$, from which it follows that $u_\lambda \in H_1$. The idea is to prove that $u_\lambda$ (or at least a subsequence thereof) converges in some sense to a certain $u$ that satisfies $Lu = -h$. In fact, in [2], it is proven that there exists such $u \in H_1$ such that $u_\lambda \to u$ strongly in $H_1$ and $Lu_\lambda \to -h$ weakly in $H_{-1}$. Since $L$ is a continuous operator, by uniqueness of limit, we have $-Lu = h$. Approximating $u$ by local functions, the lemma follows. □

4.4. *Asymmetric case for $d \geq 3$.* In dimension $d \geq 3$, a necessary and sufficient condition for a local function $v$ to be in $H_{0,-1}$ is that $\langle v \rangle_\alpha = 0$ [7]. In particular, $w_a, v_a \in H_{0,-1}$ and for $g \in \mathrm{Loc}$, $Lg \in H_{0,-1}$. As for the mean zero case, we choose $S_0 = -(L + L^*)/2$, $S_{0,N} = -(L_N + L_N^*)/2$ and apply Theorem 2. The difference here is that for $\alpha' \neq \alpha$, we have $\langle v_a \rangle_{\alpha'} \neq 0$ and so we can not invoke Theorem 3 in order to prove assumption (B). The next lemma says that condition (B) is true for this case. The proof of this lemma will be presented in the next section.

LEMMA 5. *In dimension $d \geq 3$, for a local function $h$ with $\langle h \rangle_\alpha = 0$,*

$$\lim_{N \to \infty} \left\| h - \int h \, d\mu_{N,K} \right\|_{-1,N} = \|h\|_{-1}.$$

A proof of assumption (A) for this case can be found in [6]. Once Assumptions (A) and (B) are verified, Theorem 1 follows from Theorem 2.

**5. Proof of Lemma 5.** First, note that Lemma 5 is just the generalization, in dimension $d \geq 3$, of Theorem 3 to the case in which $\langle v \rangle_\alpha = 0$ only for the fixed $\alpha \in [0,1]$. Consequently, in order to prove Lemma 5, it is enough to prove the corresponding generalizations of Lemmas 1, 2 and 3 to this case. Note that the $\|\cdot\|_{-1}$ norm depends only on the symmetric part $S$ of the generator $L$. Define the operators $S_0 = (L_0 + L_0^*)/2$ and $S_{0,N} = (L_{0,N} + L_{0,N}^*)/2$, the symmetric part of the jumps of the environment, as follows:

$$S_0 f(\xi) = \sum_{x,y \in \mathbb{Z}_*^d} s(y-x)\xi(x)(1-\xi(y))[f(\sigma^{xy}\xi) - f(\xi)],$$

$$S_{0,N} f(\xi) = \sum_{x,y \in T_{N,*}^d} s(y-x)\xi(x)(1-\xi(y))[f(\sigma^{xy}\xi) - f(\xi)].$$

The generalizations of Lemmas 1 and 3 are proven in [7] and [3].

LEMMA 6. *In dimension $d \geq 3$, if $v \in \mathrm{Loc}$ satisfies $\langle v \rangle_\alpha = 0$, then $v \in H_{0,-1}$.*

LEMMA 7. *In dimension $d \geq 3$, if $v \in \mathrm{Loc}$ and $\langle v \rangle_\alpha = 0$, then for all $\varepsilon > 0$, there exists $g_\varepsilon \in \mathrm{Loc}$ such that*

$$\|v - Sg_\varepsilon\|_{0,-1} < \varepsilon.$$

So, it only remains to prove the generalization of Lemma 2 to this case.

LEMMA 8. *Let $v$ be a local function such that $\langle v \rangle_\alpha = 0$. Define $\langle v \rangle_N = \int v \, d\mu_{N,K_N}$. In dimension $d \geq 3$,*

$$\lim_{N \to \infty} \|v - \langle v \rangle_N\|_{0,-1,N} = \|v\|_{0,-1}.$$

PROOF. Using the variational formula for $\|v\|_{0,-1}$, it is not hard to prove that

$$\liminf_{N \to \infty} \|v - \langle v \rangle_N\|_{0,-1,N} \geq \|v\|_{0,-1}.$$

In fact, by definition, for all $\varepsilon > 0$, there exists a local function $f_\varepsilon$ such that

$$\begin{aligned}
\|v\|_{0,-1}^2 &\leq 2\langle v, f_\varepsilon \rangle_\alpha - \langle f_\varepsilon, -S_0 f_\varepsilon \rangle_\alpha + \varepsilon \\
&= \lim_{N \to \infty} \{2\langle v - \langle v \rangle_N, f_\varepsilon \rangle_N - \langle f_\varepsilon, -S_{0,N} f_\varepsilon \rangle_N\} + \varepsilon \\
&\leq \liminf_{N \to \infty} \sup_f \{2\langle v - \langle v \rangle_N, f \rangle_N - \langle f, -S_{0,N} f \rangle_N\} + \varepsilon \\
&= \liminf_{N \to \infty} \|v - \langle v \rangle_N\|_{0,-1,N}^2 + \varepsilon.
\end{aligned}$$

The converse inequality is harder to prove. The idea is to approximate $v$ in $H_{0,-1}$ by local functions with mean zero for all densities $\alpha \in [0,1]$. The proof requires two auxiliary lemmas. The first is just a version of Lemma 3.6 of [3].

LEMMA 9. *Let $w$ be a local function with $\langle w \rangle_\alpha = 0$ for all $\alpha \in [0,1]$. Let $\{f_N\}_N$ be a sequence of functions defined in $H_{0,1,N}$ such that*

$$\langle f_N, -S_{0,N} f_N \rangle_N \leq 1,$$
$$\lim_{N \to \infty} \langle w, f_N \rangle_N = A.$$

*Then there exist $f \in H_{0,1}$ and subsequence $N'$ such that $\langle w, f \rangle_\alpha = A$, $\langle f, -S_0 f \rangle_\alpha \leq 1$ and for all local functions $h$ with $\langle h \rangle_\alpha = 0$, for each $\alpha \in [0,1]$,*

$$\lim_{N' \to \infty} \langle f_{N'}, h \rangle_N = \langle f, h \rangle_\alpha.$$

Before stating the second auxiliary lemma, we need to introduce some notation. Let $\Lambda_N = \{-N+1, \ldots, N\}^d \setminus \{0\}$ be the cube of radius $N$. Note that $\Lambda_N \neq T_{N,*}^d$ because $\Lambda_N$ has no periodic boundary conditions. For each $x \in \mathbb{Z}_*^d$, define $\theta_x(\xi) =: \xi(x)$ and for each $l > 0$, define $\varphi_l(\xi) = \sum_{x \in \Lambda_l} \xi(x)$. Let $\mathcal{F}_{\Lambda_N}$ be the $\sigma$-algebra generated by $\varphi_l$ and $\{\theta_x ; x \in \Lambda_N^c\}$. For $l > 0$ such that $\mathrm{supp}(v) \subseteq \Lambda_l$, define $v_l = \mathbb{E}[v|\mathcal{F}_{\Lambda_l}]$. Note that there is a natural way to define $v_l$ that does not depend on the particular value of $\alpha$. The following lemma is an easy consequence of the equivalence of ensembles:

LEMMA 10. *Fix positive integers $l, q$ such that $\mathrm{supp}(v) \subseteq \Lambda_l$ and $q > 2$. Define $g_n = v_{lq^n}$. Then there is a finite constant $\kappa$ such that:*

(i) $\langle (g_n - g_{n-1})^2 \rangle_\alpha \leq \kappa (lq^n)^{-d}$,
(ii) $\langle (g_n - g_{n-1})^2 \rangle_N \leq \kappa (lq^n)^{-d}$.

The proof is as follows. For each $N$, there exists a function $f_N \in H_{1,N}$ such that

$$\langle f_N, -S_{0,N} f_N \rangle_N \leq 1 \quad \text{and} \quad \|v - \langle v \rangle_N\|_{0,-1,N} = \langle f_N, v - \langle v \rangle_N \rangle_N.$$

Consider a subsequence $\tilde{N}$ such that

$$\lim_{\tilde{N} \to \infty} \|v - \langle v \rangle_{\tilde{N}}\|_{0,-1,\tilde{N}} = \limsup_{N \to \infty} \|v - \langle v \rangle_N\|_{0,-1,N} =: A.$$

By Lemma 9, there exists a function $f \in H_1$ and a subsubsequence $N'$ such that $\langle f_{N'}, h \rangle_{N'} \to \langle f, h \rangle_\alpha$ for all local functions $h$ with mean zero for each $\mu_\alpha$. In particular,

$$\lim_{N' \to \infty} \langle f_{N'}, v - v_l \rangle_{N'} = \langle f, v - v_l \rangle_\alpha.$$

Let $l, q > 2$ be fixed. Define, as in Lemma 9, $g_n = v_{lq^n}$. In order to make notation simpler, suppose that $N' = lq^n$ and denote $N'$ simply by $N$. If this is not the case, then the required changes are straightforward. We then have that

$$\langle f_N, v - \langle v \rangle_N \rangle_N = \langle f_N, v - v_l \rangle_N + \langle f_N, v_l - \langle v \rangle_N \rangle_N$$
$$= \sum_{k=1}^{n} \langle f_N, g_{k-1} - g_k \rangle_N + \langle f_N, v - v_l \rangle_N.$$

Define $\mathcal{L}_k$ as the generator of an exclusion process in $\Lambda_{lq^k}$. Note that due to the boundary effects, $\mathcal{L}_{lq^k} \neq S_{0,lq^k}$. We see that $\langle v - v_l \rangle_\alpha = 0$, $\langle g_{k-1} - g_k \rangle_\alpha = 0$ for all $\alpha \in [0,1]$. By linear algebra, there exists a local function $G_k$ defined in $\{0,1\}^{\Lambda_{lq^k}}$ such that $g_{k-1} - g_k = \mathcal{L}_k G_k$. Therefore,

$$\sum_{k=1}^{n} \langle f_N, g_{k-1} - g_k \rangle_N = \sum_{k=1}^{n} \langle f_N, \mathcal{L}_k G_k \rangle_N$$
$$= \sum_{k=1}^{n} \sum_{b \in \Gamma_k} \langle \nabla_b f_N, \nabla_b G_k \rangle_N,$$

where $\sum_{b \in \Gamma_k}$ indicates a sum over all bonds $b = \langle xy \rangle$ such that $x, y \in \Lambda_{lq^k}$ and $\nabla_b g = s(y - x)^{1/2} [g(\sigma^{xy} \eta) - g(\eta)]$.

Choose $a_k = \varepsilon 2^k$. By Cauchy's inequality with weights $a_k$, we have

$$\left| \sum_{k=1}^{n} \langle f_N, g_{k-1} - g_k \rangle_N \right|$$
$$\leq \sum_{k=1}^{n} \sum_{b \in \Gamma_k} \frac{1}{a_k} \langle (\nabla_b f_N)^2 \rangle_N + a_k \langle (\nabla_b G_k)^2 \rangle_N$$
$$\leq \sum_{b \in \Gamma_n} \sum_{k : b \in \Gamma_k} \frac{1}{a_k} \langle (\nabla_b f_N)^2 \rangle_N + \sum_{k=1}^{n} \sum_{b \in \Gamma_k} a_k \langle (\nabla_b G_k)^2 \rangle_N$$
$$\leq \frac{1}{\varepsilon} \sum_{b \in \Gamma_n} \langle (\nabla_b f_N)^2 \rangle_N + \sum_{k=1}^{n} a_k \langle g_k - g_{k-1}, -\mathcal{L}_k^{-1}(g_k - g_{k-1}) \rangle_N$$
$$\leq \frac{1}{\varepsilon} \langle f_N, -\mathcal{L}_n f_N \rangle_N + \varepsilon \sum_{k=1}^{n} 2^k \langle g_k - g_{k-1}, -\mathcal{L}_k^{-1}(g_k - g_{k-1}) \rangle_N$$
$$\leq \frac{1}{\varepsilon} \langle f_N, -S_{0,lq^n} f_N \rangle_N + \varepsilon \sum_{k=1}^{n} 2^k C \cdot 2^k (lq^k)^2 \langle (g_{k-1} - g_k)^2 \rangle_N,$$

where, in the last line, we have used the spectral gap inequality for the exclusion process [4].

Using Lemma 10 and minimizing in $\varepsilon$, we have

$$\left|\sum_{k=1}^n \langle f_N, g_{k-1} - g_k\rangle_N\right| \leq \frac{1}{\varepsilon} + \varepsilon \sum_{k=1}^n C\kappa \cdot 2^k (lq^k)^{2-d}$$

$$\leq \frac{1}{\varepsilon} + \varepsilon\left[\frac{C\kappa l^{2-d}}{1 - 2q^{2-d}}\right]$$

$$\leq 2\sqrt{\frac{C\kappa l^{2-d}}{1 - 2q^{2-d}}} \leq C_1 l^{(2-d)/2}.$$

By the law of large numbers, as $l \to \infty$, $v_l \to 0$ $\mu_\alpha$-a.s. and in $\mathcal{L}^2(\mu_\alpha)$. We also have that

$$\|g_k - g_{k-1}\|_{0,-1}^2 = \langle g_k - g_{k-1}, (-S_0)^{-1}(g_k - g_{k-1})\rangle_\alpha$$

$$\leq \langle g_k - g_{k-1}, (-\mathcal{L}_{k+1})^{-1}(g_k - g_{k-1})\rangle_\alpha$$

$$\leq C(lq^{k+1})^2 \langle (g_k - g_{k-1})^2\rangle_\alpha$$

$$\leq C\kappa q^2 (lq^k)^{2-d}.$$

Therefore, the sequence $\{g_k - g_{k-1}\}_k$ is absolutely summable and there exists $g \in H_{0,-1}$ such that

$$\lim_{n\to\infty}(v_l - v_{lq^n}) = \sum_{k=1}^\infty g_k - g_{k-1} = g.$$

On the other hand, we know that $v_{lq^n} \to 0$ in $\mathcal{L}^2(\mu_\alpha)$, from which it follows that $\langle F, v_{lq^n}\rangle_\alpha$ goes to zero for all $F \in \mathcal{L}^2(\mu_\alpha)$ and $v_{lq^n} \to v_l - g$ in $H_{0,-1}$. From this, $\langle F, v_{lq^n}\rangle_\alpha \to \langle F, v_l - g\rangle_\alpha$ for all $F \in H_{0,1}$. Since $D(S_0) \subseteq \mathcal{L}^2(\mu_\alpha) \cap H_{0,-1}$ and $D(S_0)$ is dense in $H_{0,-1}$, we have $g = v_l$.

As before, by using part (i) of Lemma 10, we can prove that there exists a constant $C_2$ such that

$$|\langle f, v_l\rangle_\alpha| \leq C_2 \cdot l^{(2-d)/2}.$$

Combining both inequalities, we see that

$$\limsup_{N\to\infty} \|v - \langle v\rangle_N\|_{0,-1,N} = \limsup_{N\to\infty} \langle f_N, v - \langle v\rangle_N\rangle_N$$

$$= \limsup_{N\to\infty}\{\langle f_N, v - v_l\rangle_N + \langle f_N, v_l - \langle v\rangle_N\rangle_N\}$$

$$\leq \langle f, v - v_l\rangle_\alpha + C_1 \cdot l^{(2-d)/2}$$

$$\leq (C_1 + C_2)l^{(2-d)/2} + \langle f, v\rangle_\alpha.$$

Since $d \geq 3$ and $l$ is arbitrary, we have

$$\limsup_{N\to\infty} \|v - \langle v\rangle_N\|_{0,-1,N} \leq \langle f, v\rangle_\alpha \leq \|f\|_{0,1} \cdot \|v\|_{0,-1} \leq \|v\|_{0,-1}. \qquad \square$$

**Acknowledgments.** The author would like to thank his thesis advisor, Claudio Landim, for suggesting this problem and for his valuable support during the undertaking. Thanks are due to the referee for valuable suggestions that lead to this final version of the article.

## REFERENCES

[1] KIPNIS, C. and VARADHAN, S. R. S. (1986). Central limit theorem for additive functionals of reversible Markov process and applications to simple exclusion process. *Comm. Math. Phys.* **106** 1–19. MR0834478

[2] LANDIM, C., OLLA, S. and VARADHAN, S. R. S. (2000). Asymptotic behavior of a tagged particle in simple exclusion process. *Bol. Soc. Brasil. Mat.* **31** 241–275. MR1817088

[3] LANDIM, C., OLLA, S. and VARADHAN, S. R. S. (2002). Finite-dimensional approximation of the self-diffusion coefficient for the exclusion process. *Ann. Probab.* **30** 483–508. MR1905849

[4] QUASTEL, J. (1992). Diffusion of color in the simple exclusion process. *Comm. Pure Appl. Math.* **40** 623–679. MR1162368

[5] SAADA, E. (1987). A limit theorem for the position of a tagged particle in a simple exclusion process. *Ann. Probab.* **15** 375–381. MR0877609

[6] SETHURAMAN, S., VARADHAN, S. R. S. and YAU, H.-T. (2000). Diffusive limit of a tagged particle in asymmetric simple exclusion process. *Comm. Pure Appl. Math.* **53** 972–1006. MR1755948

[7] SETHURAMAN, S. and XU, L. (1996). A central limit theorem for reversible exclusion and zero-range particle systems. *Ann. Probab.* **24** 1842–1870. MR1415231

[8] VARADHAN, S. R. S. (1995). Self-diffusion of a tagged particle in equilibrium for asymmetric mean zero random walks with simple exclusion. *Ann. Inst. H. Poincaré Probab. Statist.* **31** 273–285. MR1340041

IMPA
ESTRADA DONA CASTORINA 110
CEP 22460 RIO DE JANEIRO
BRAZIL
E-MAIL: monets@impa.br